\documentclass{article}

\def \Z {{\mathbb {Z}}}

\def \P {\mathbf {P}}

\def \R {{\mathbb {R}}}

\usepackage[tbtags]{amsmath}
\usepackage{amsfonts,amssymb}

\title{Explicit Infinite Mixing Automorphisms \\
       with Simple Spectra of Symmetric  Squares}
\author{ Sofia V. Vereshchagina, Valery V. Ryzhikov }
\date{}

\begin{document}

\date{}

\maketitle

\Large
\begin{abstract}

We present mixing ergodic automorphisms of a space with sigma-finite measure whose symmetric tensor square has single spectrum. This property is of interest in connection with problems of A.N. Kolmogorov and V.A. Rokhlin on the spectra of ergodic automorphisms. Explicit examples are given in the class of rank-one constructions.

\vspace{2mm}
\end{abstract}

\vspace{5mm}

\section{Introduction}
In the 1960s, A.N. Kolmogorov conjectured that the spectral measure of an ergodic automorphism subordinates its convolution.
This phenomenon is called the group property of the spectrum.

Kolmogorov's conjecture was confirmed for automorphisms with discrete spectrum, Gaussian systems, K-automorphisms, and dynamical systems proposed by Ya.G. Sinai in \cite{Si}.

In the class of non-mixing automorphisms, counterexamples were given in the works of A.B. Katok, A.M. Stepin \cite{KS}, and V.I. Oseledets \cite{O}. For weakly mixing automorphisms, a negative answer was given thanks to the solution of a problem known in ergodic theory as V.A. Rokhlin's problem on the homogeneous spectrum of an ergodic automorphism. This topic is covered, for example, in D. V. Anosov's lectures \cite{An}. In \cite{07}, the existence of a mixing automorphism $T$ of a probability space was shown for which the product $T\times T$ has a homogeneous spectrum of multiplicity 2. A consequence of this fact is the absence of a group property for such   $T$ (otherwise, the spectrum multiplicity of the product $T\times T$ would be greater than 2). Thus, in the class of mixing automorphisms, Kolmogorov's problem was solved along with Rokhlin's problem. The paper \cite{07} did not provide specific examples. Explicit constructions are described in \cite{24}. These include staircase constructions with sufficiently slow power-law growth of the parameters $r_j$. The mixing property for such constructions was proved by T. Adams \cite{Ad}.
We note that proving the mixing property for staircase  constructions is nontrivial. In our case, mixing follows clearly from the definition of the constructions themselves. We slightly modify the so-called Sidon transformations to obtain the desired spectral property.

The automorphism $T$ and the change-of-variables operator induced by it in $L_2(\mu)$
are denoted identically. In the case of automorphisms of a probability space, the mixing property is equivalent to $\tilde T^n\to_w 0$, where $\tilde T$ is the restriction of $T$ to a space orthogonal to the constants. In our case, nonzero constants do not lie in $L_2$, so $T^n\to_w 0$, $n\to\infty$, is the definition of the mixing property.
The ordinary shift by $\R$ is a mixing transformation. Here, of course, the physical analogy is lost, but it is restored by considering Poisson
superstructures over such formally mixing transformations (see \cite{26}).

Let us briefly outline the method we use to obtain the simple spectrum of a symmetric degree of order 2.

Let $U:H\rightarrow H$ be a unitary operator on a separable Hilbert space $H$. A cyclic subspace of $U$ with cyclic vector $h\in H$ is defined as
the subspace $C_h:=\overline{<U^nh\:|\:n\in\mathbb{Z}>}$ (the closure of the linear span of the vectors
$U^nh$). If $H=C_h$ for some $h$, $U$ is said to have simple (single) spectrum. We consider the operator $T$ defined below on $L_2(\mu)$ with a cyclic vector $f$ and prove that the product $h=f\otimes f$ is a cyclic vector for the symmetric tensor
square $U=T\odot T$ acting on the space $H=L_2(\mu)\odot L_2(\mu)$. The operator $T$ is induced by an automorphism $T$ of the space $(X,\mu)$ with sigma-finite measure $\mu$. A modified Sidon automorphism is used as a suitable $T$.
It is chosen such that
the cyclic space $C_{f\otimes f}$ contains all vectors of the form $F_n= T^n f\otimes f +f\otimes T^nf$, $n=0, 1, 2,\dots $. This membership is established using approximation considerations. We will show that specially chosen linear combinations of vectors $U^ih$ can approximate vectors $F_n$ arbitrarily well.
As a result, we find that all vectors of the form $U^mF_n$ belong to the space $C_{f\otimes f}$.
The closure of linear combinations of vectors $U^mF_n$ ($m\in \Z$, $n=0,1,2,\dots$) is the space $L_2(\mu)\odot L_2(\mu)$, which, by virtue of the above, coincides with $C_{f\otimes f}$. This establishes the desired simple spectrum of the operator $U=T\odot T$.

\section{Automorphisms of the class ${\cal T}_2$}
Let us recall the definition of a rank-one automorphism (we will find the examples we need among them).
We fix a natural number $h_1\geq 4$ (the height of the tower at stage $j=1$), a sequence $r_j$ (the parameter $r_j$ is the number of columns into which the tower of stage $j$ is virtually cut), and a sequence of integer vectors (superstructure parameters)
$$ \bar s_j=(s_j(1), s_j(2),\dots, s_j(r_j-1),s_j(r_j)).$$
Below is a description of the construction of a measure-preserving transformation, which is completely determined by the parameters $h_1$, $r_j$, and $\bar s_j$.

At stage $j=1$, the half-interval $E_1$ is defined. At step $j$, a system of disjoint half-intervals is defined:
$$E_j, TE_j, T^2E_j,\dots, T^{h_j-1}E_j,$$
where the transformation $T$ on $E_j, TE_j,\dots, T^{h_j-2}E_j$
is a parallel translation. This set of half-intervals is called the tower of step $j$; their union is denoted by $X_j$ and is also called a tower. A half-interval included in the tower of step $j$ is also called a tower floor of step $j$, or simply a floor when the step number is not specified.

Let's represent $E_j$ as a disjoint union of $r_j$ intervals
$$E_j^1,E_j^2,E_j^3,\dots, E_j^{r_j}$$ of the same measure (length).
For each $i=1,2,\dots, r_j$, we consider the so-called column
$$E_j^i, TE_j^i ,T^2 E_j^i,\dots, T^{h_j-1}E_j^i.$$
We denote the union of these intervals by $X_{i,j}$; we will also call it a column.

To each column $X_{i,j}$ we add $s_j(i)$ disjoint half-intervals of the same measure as $E_j^i$, obtaining a set
$$E_j^i, TE_j^i, T^2 E_j^i,\dots, T^{h_j-1}E_j^i, T^{h_j}E_j^i, T^{h_j+1}E_j^i, \dots, T^{h_j+s_j(i)-1}E_j^i$$
(all these sets are disjoint).
Denoting $E_{j+1}= E^1_j$, for $i<r_j$ we set
$$T^{h_j+s_j(i)}E_j^i = E_j^{i+1}.$$
From this point on, the set of superstructured columns is considered as a tower of stage $j+1$,
consisting of the half-intervals
$$E_{j+1}, TE_{j+1}, T^2 E_{j+1},\dots, T^{h_{j+1}-1}E_{j+1},$$
where
$$ h_{j+1} =h_jr_j +\sum_{i=1}^{r_j}s_j(i).$$

The definition of the transformation $T$ at stage $j$ is preserved in the following stages. As a result, we obtain the space $X=\cup_j X_j$ and an invertible transformation $T:X\to X$ that preserves the standard Lebesgue measure on $X$.
If
$$\sum\limits_{j} {(s_j(1)+s_j(2)+\dots+s_j(r_j))}/{h_jr_j }=\infty,$$
the measure of $X$ is infinite, this is the case we will consider below.

{ \bf $k$-Sidon automorphisms.}
Let a rank-one construction $T$ have the following property: \it the intersection
$X_j\cap T^mX_j$ for $h_{j}<m\leq h_{j+1}$ can be contained
in the union of at most $k$ columns $X_{i,j}$ of the tower $X_j$. \rm Such a construction is called \it $k$-Sidon. \rm
Note that for such $T$ and a set $A$ consisting of the tower floors of stage $j_0$, $$\mu(A\cap T^mB)\leq k\mu(A)/r_j$$ holds for all $m\in [h_j, h_{j+1}],$ $j>j_0$.
It follows that the $k$-Sidon construction is mixing, which by definition means
$$\mu(A\cap T^mB)\to 0, m\to\infty$$ for all sets $A,B$ of finite measure.

\bf Automorphisms of the class ${\cal T}_2$. \rm
Let the parameters $h_1,\:r_j,\:\overline{s_j}$ be given.
For natural numbers $r,n$, define the sets
$$J(r,n):=\{j\:|r_j=r,s_{j}(r_j-1)=n\}.$$
We say that an automorphism $T$ of rank 1 belongs to the class ${\cal T}_2$ if the parameters satisfy the following conditions:

\vspace{2mm}
$r_j$  is nondecreasing and takes all values starting with $r_1>4$;

$J(r,n)= \varnothing$ for $r<n$; \ \ \ $\overline{\lim\limits_{r\to\infty}}\frac{|J(r,n)|}{r^4}=+\infty$ for all $n$

$s_j(1)\geq 4h_j$;  
$s_j(i)\geq 4s_j(i-1)$ for $i=2,...,r_j-3$;

$s_j(r_{j}-2)=0$;\ \ \  $s_j(r_j)\geq 4s_j(r_j-3)$.

\vspace{2mm}
Under the stated conditions, the phase space measure of the corresponding construction is infinite. Such constructions are 2-Sidon and are mixing. Let us formulate the main result.

\vspace{3mm}
\textbf{Theorem.} \it Automorphisms $T$ of class $ {\cal T}_2$ are mixing and their spectra do not have Kolmogorov's group property. Spectra  of their symmetric tensor squares $T\odot T$ are simple, so the products $T\otimes T$ have  homogeneous spectra of multiplicity 2. \rm

\newpage 
\section{Proof }
Let $T\in {\cal T}_2$, and fix some floor $A$ of tower $X_1$. We set $A_{k,j}:=A\cap X_{k,j},$ $k=1,\dots,r_j$. We will also denote by $A$ and $A_{k,j}$ the characteristic functions of the corresponding sets.

\vspace{3mm}
\textbf{Lemma 1.} \it Let $j\in J(r,n).$ For $m=0$ or $n$, we have
$$\big(T^{h_j}A,T^{-m}A\big)=\frac{1}{r_j}\big(A+T^{-s_j(r_j-1)}A,T^{-m}A\big)$$ \rm

\vspace{6mm}
The proof of this and the following lemmas is verified directly;
we omit this verification in this text.

\vspace{6mm}
\textbf{Lemma 2.} For any distinct $i, j \in J(r, n)$, the following equality holds:
\[
(T^{h_j}A, T^{h_i}A) = \frac{1}{r^2} \big( (I + T^{-n})A, (I + T^{-n})A \big).
\]

\vspace{6mm}
\textbf{Lemma 3.} \it $$\left\| \frac{1}{|J(r,n)|}\sum\limits_{j\in J(r,n)}T^{h_j}A\otimes T^{h_j}A - \frac{1}{r^2}(I+T^{-n})A\otimes(I+T^{-n})A\right\|^2=$$
$$=\frac{1}{|J(r,n)|}\Big(\big( A,A\big)^2 - \frac{1}{r^4}\big(A+T^{-n}A,A+T^{-n}A \big)^2\Big).$$ \rm

\newpage
{\bf Proof of  Theorem.}
We need to show that for an automorphism $T$ of class $ {\cal T}_2$, the spectrum of the product $T\odot T$ is single-valued
To do this, it is sufficient to verify that the functions $F_n=(I+T^{-n})f\otimes(I+T^{-n})f$ can be approximated arbitrarily accurately by vectors from $C_{f\otimes f}$.

For integers $r,n$ such that the set $J(r,n)$ is nonempty, we set $${\P_{r,n}}:=\frac{1}{|J(r,n)|}\sum\limits_{j\in J(r,n)}T^{h_j}\otimes T^{h_j}.$$

By Lemma 3, we have the inequality $$\left\|r^2{\P_{r,n}}f\otimes f - (I+T^{-n})f\otimes(I+T^{-n})f\right\|^2\leq \frac{Cr^4}{|J(r,n)|}$$ for some positive constant $C.$

By the definition of the class ${\cal T}_2$, for every $n$ from the sequence $\{r_j\}_{j=1}^{\infty}$ one can choose a subsequence $\{r(j,n)\}_{j=1}^{\infty}$ such that the condition $$\lim\limits_{j\to \infty}\frac{{r(j,n)}^4}{|J(r(j,n),n)|}=0.$$ Then we obtain 
$$\left\|{r(j,n)}^2{\P_{r(j,n),n}}f\otimes f - (I+T^{-n})f\otimes(I+T^{-n})f\right\|=o(1),\:j\to\infty,$$
this completes the proof of the theorem. 

\vspace{3mm}
\bf Remark. \rm In this note, we confine ourselves to studying the spectral properties of constructions of the class ${\cal T}_2$.
Using $k$-Sidon constructions, one can construct explicit examples such that both the spectal simplicity  for symmetric tensor powers of small order and the absolute continuity of  spectra for large order powers  are realized. A particular result in this direction is formulated in \cite{25}.

\end{document}